\documentclass[12pt]{amsart}

\usepackage{amsmath,xspace,amssymb,mathrsfs}
\usepackage{color}

\input xy
\xyoption{all}
\xyoption{2cell}
\UseAllTwocells
\CompileMatrices

\newcommand{\mr}{\operatorname{\mathfrak{mr}}}
\newcommand{\Cl}{\operatorname{Clop}}
\newcommand{\Spec}{\operatorname{Spec}}
\renewcommand{\phi}{\varphi}
\newcommand{\Sp}{\operatorname{Sp}}
\newcommand{\Gal}{\operatorname{Gal}}
\newcommand{\hm}{\operatorname{Hom}}

  \newtheorem{proposition}{Proposition}[section]
  \newtheorem{lemma}[proposition]{Lemma}
  
  \newtheorem{corollary}[proposition]{Corollary}
  \newtheorem{theorem}[proposition]{Theorem}

  \theoremstyle{definition}
  \newtheorem{definition}[proposition]{Definition}

  \newtheorem{remark}[proposition]{Remark}

\begin{document}

\title{On the category of profinite spaces as a reflective subcategory}

\author{Abolfazl Tarizadeh}
\address{ Department of mathematics, Institute for  Advanced Studies in Basic Sciences(IASBS),
 P. O. Box 45195-1159, Zanjan, Iran.}
\email{abtari@iasbs.ac.ir}

\date{}

\begin{abstract} In this paper by using the ring of real-valued continuous functions $C(X)$, we prove a theorem in profinite spaces which states that for a compact Hausdorff space $X$, the set of its connected components $X/_{\sim}$ endowed with some topology $\mathscr{T}$ is a profinite space. Then we apply this result to give an alternative proof to the fact that the category of profinite spaces is a reflective subcategory in the category of compact  Hausdorff spaces. Finally, under some circumstances on a space $X$, we compute the connected components of the space $t(X)$ in terms of the ones of the space $X$.

\end{abstract}

\maketitle

\section{Introduction}

A profinite space is a compact Hausdorff
and totally disconnected topological space. In other words, a space $X$ is profinite if there exists an inverse system of finite discrete spaces for which its inverse limit is homeomorphic to $X$, consider \cite[Section 3.4]{Borceux Janelidze}. Recall that a profinite group is a topological group whose underlying space is a profinite space. \\

There are  interesting examples of profinite spaces and profinite groups which arise from algebraic geometry, Galois theory and topology. For instance, for any field $K$ its absolute Galois group $\Gal(K^{s}/K)$ is a profinite group, or more generally the \'{e}tale fundamental group $\pi_{1}(X,\overline{s})$ of a connected scheme $X$ on a geometric point $\overline{s}:\Spec(\Omega)\rightarrow X$ is a profinite group \cite[Theorem 5.4.2]{Szamuely}. \\

 Stone's duality says us that any profinite space $X$ is of the form $X=\Spec(B)$ for some Boolean algebra $B$ (\cite[Theore 4.1.16]{Borceux Janelidze}). \\
  There is also another characterization of profinite spaces due to Craven \cite{Craven}, where he proves that each profinite space is homeomorphic to the space $X(F)$ for some formally real field $F$ and $X(F)$ denotes the set of orderings of the field $F$ endowed with some topology.  \\

In partial of  this paper, we give an alternative proof to the fact that the category of profinite spaces is a reflective subcategory in category of compact Hausdorff spaces. We prove this by using spectra of the Boolean algebras and rings of real-valued continuous functions $C(X)$.\\Thanks to the Stone-\v{C}ech compactification functor, one can prove that the category of compact Hausdorff spaces is a reflective subcategory in the category of Tychonoff spaces. Moreover,
the category of Tychonoff spaces itself is reflective in the category of topological spaces. \\

  Section 2, contains some preliminaries which will be required in Section 3. In  Section 3, we will use spectra of the Boolean algebras to compute  the connected components of the spectra of a commutative ring in terms of the its max-regular ideals (Theorem \ref{max-reg}). As an application of this result and also by using some properties of the rings of real-valued continuous functions $C(X)$, we give an alternative proof to a theorem in profinite spaces which states that for
a compact  Hausdorff space $X$, the set of its connected components $X/_{\sim}$ endowed with some topology $\mathscr{T}$ is a profinite space (Theorem \ref{cfunction}). This theorem leads us to construct a covariant functor from the category of compact Hausdorff spaces $\mathscr{K}$ to the category of profinite spaces $\mathscr{P}$, then we will use this categorical construction to show that  the category $\mathscr{P}$ is a reflective subcategory in the category $\mathscr{K}$ (Theorem \ref{alternative}). Finally, under some circumstances on a space $X$, we compute the connected components of $t(X)$ in terms of the connected components of $X$(Theorem \ref{conn component}). Consider Section 5.7 of the book \cite{Borceux Janelidze} to see another proof of the theorems \ref{cfunction} and  \ref{alternative}. In that book, These theorems are proven by means of the nearness relation. \\

\section{ Preliminaries }
 In this Section, for convenience of the reader and for the sake of completeness we collect some preliminaries which will be required in the next section. For more details on the spectra of the Boolean algebras we reference the reader to the book \cite[Sections 4.1 and 4.2] {Borceux Janelidze}. \\

\begin{definition}  A topological space $X$ is said to be a profinite space if it is compact Hausdorff and totally disconnected. Totally disconnectedness means that there is no connected subset in $X$ other than its points.\\
\end{definition}

\begin{definition} A Boolean algebra is a structure $(B,\vee ,\wedge ,c,0,1)$ with two binary operations $\vee,\wedge: B\times B\rightarrow B$, a unary operation $c:B\rightarrow B$ and two distinguished elements 0 and 1 in $B$, such that for all $x,y$ and $z$ in $B$,\\
(i) the  binary operations $\vee$ and $\wedge$ are commutative and associative,\\
(ii) $x\wedge(y\vee z)=(x\wedge y)\vee(x\wedge z)$, \\
(iii) $x\vee(y\wedge z)=(x\vee y)\wedge(x\vee z)$,\\
(iv) $x\vee(x\wedge y)=x$, $x\wedge (x\vee y)=x$, \\
(v) $x\vee x^{c}=1$, $x\wedge x^{c}=0$. \\
\end{definition}

 For the sake of simplicity, the Boolean algebra $(B,\vee ,\wedge ,c,0,1)$ is only denoted  by $B$. The relation $x\leq y \Longleftrightarrow x\wedge y=x \Longleftrightarrow x\vee y=y$ puts a partial ordering on $B$. \\

  A morphism between the Boolean algebras is a map $\varphi:B\rightarrow B'$ which preserves the binary and unary operations.

\begin{definition} Let $B$ be a Boolean algebra. A filter in $B$ is a subset $F\subseteq B$ such that,\\
(i) $1\in F$, \\
(ii) if $x,y\in F$ then $x\wedge y\in F$, \\
(iii) if $x\in F$ and $x\leq y$ then $y\in F$. \\

The filter $F$ is called proper if $0\notin F$. Each maximal element of the poset of proper filters ordered by inclusion is called an ultrafilter.  \\
\end{definition}

\begin{definition} \label{stone basis} For given Boolean algebra $B$, we denote by $\Spec(B)$ the set of all its ultrafilters. For every filter $H$ on $B$ consider $\mathscr{O}_{H}=\{F \in \Spec(B)$ $|$ $H\nsubseteq F\}.$ The subsets $\mathscr{O}_{H}$ as open subsets constitute a topology on $\Spec(B)$. This topology is usually called  the Stone topology  and  the space $\Spec(B)$ is called the spectrum of $B$. Also the collection of all $\mathscr{O}_{b}=\{F \in \Spec(B)$ $|$ $ b\notin F\}$ where $b\in B$, as open subsets constitute a basis for the Stone topology. \\

\end{definition}

\begin{lemma}\label{bool-profinite} The spectrum of a Boolean algebra is a profinite space.
\end{lemma}

  {\bf Proof.} For the proof see \cite[Proposition 4.1.11] {Borceux Janelidze}.  $\Box$  \\

All rings that we consider in this paper are commutative with the identity. \\

\begin{lemma} For a ring $R$, then the set of its idempotents $I(R)$ with the operations $e\vee e'=e+e'-ee'$, $e\wedge e'=ee'$ and $e^{c}=1-e$ constitute a Boolean algebra.
\end{lemma}
{\bf Proof.} For the proof see \cite[Proposition 4.2.1] {Borceux Janelidze}.  $\Box$  \\

\begin{remark} \label{sp} We denote by $\Sp(R)$ the spectrum of the Boolean algebra $I(R)$ which is a profinite space according to the  lemma \ref{bool-profinite}.
\end{remark}

\begin{definition} In a ring $R$, an ideal $I$ of $R$ is called a regular ideal if it is generated by a set of idempotents of $R$. If moreover $I\neq R$, then we call it a proper regular ideal.
\end{definition}

\begin{lemma} \label{regularidem} For any ideal $I$ of  $R$, then $I$ is a regular ideal if and only if for any $a\in I$ there exists an idempotent $e\in I$ such that $a=ea$.
\end{lemma}

{\bf Proof.} See \cite[Lemma 4.2.7]{Borceux Janelidze} for its proof. $\Box$ \\

\begin{definition}\label{maxreg def} For a ring $R$, Any maximal element of the poset of proper regular ideals of $R$ ordered by inclusion is called a max-regular ideal of $R$. We denote by $\mr(R)$ the set of all its max-regular ideals. \\ The subsets $\mathscr{O}_{I}=\{M\in \mr(R)$ $|$ $I\nsubseteq M\}$ where $I$ is a regular ideal of $R$, as open subsets constitute a topology on $\mr(R)$.
\end{definition}

\begin{lemma} \label{mrprofinite} For a ring $R$, then there exists a natural map $\phi: \Sp(R) \rightarrow \mr(R)$ which is a homeomorphism.
\end{lemma}

{\bf Proof.} Consider \cite[Corollary 4.2.11]{Borceux Janelidze} for the proof. $\Box$  \\

\begin{remark} \label{fineremark} In the above lemma, the explicit description of the map $\phi$ is as follows. Each element $F \in \Sp(R)$  is an ultarfilter of the Boolean algebra $I(R)$,  define $\phi(F)=M$ where $M=\langle e\in I(R)$ $|$ $1-e\in F \rangle$. \\In particular, as a consequence of the above lemma, the subsets $\mathscr{O}_{e}=\{M\in \mr(R)$ $|$ $e\notin M\}$ where $e\in R$ is an idempotent, constitute a basis for the topology of $\mr(R)$.
\end{remark}

\begin{definition} For any ring $R$, set $\Spec(R)=\{\mathfrak{p}\subset R$ $|$ $\mathfrak{p}$ is a prime ideal of $R\}$. Then the subsets $\mathscr{O}_{I}=\{\mathfrak{p}\in \Spec(R)$ $|$  $I\nsubseteq \mathfrak{p}\}$ where $I$ is an ideal of $R$, constitute a topology for $\Spec(R)$. This topology  is called the Zariski topology. We denote by $V(I)$ the complement of $\mathscr{O}_{I}$. Also, the subsets $D(f)=\{\mathfrak{p}\in\Spec(R)$ $|$ $f\notin \mathfrak{p}\}$ where $f\in R$, as open subsets constitute a basis for the topology.
\end{definition}

For given topological space $X$, denote by $\Cl(X)$ the set of all subsets of $X$ which are both open and closed in $X$. Each element of $\Cl(X)$ is called a clopen of $X$.\\

Throughout this paper, for any topological space $X$, we denote by $X/_{\sim}$ the set of all its connected components. In the next section, the set $X/_{\sim}$ will be equipped with some topology $\mathscr{T}$ which is coarser than the quotient topology and so the canonical projection $\pi:X\rightarrow (X/_{\sim}, \mathscr{T})$ will remain continuous for this new topology. Recall that for each $x\in X$, $\pi(x)$ is defined to be the connected component of $X$ which contains $x$.  \\

For any two covariant functors $F:\mathscr{C}\rightarrow\mathscr{D}$ and $G:\mathscr{D}\rightarrow\mathscr{C}$, the covariant functor $\hm_{\mathscr{D}}(F(-),-):\mathscr{C}^{op}\times\mathscr{D}\rightarrow\mathbf{Set}$ is called the left hom-set adjunction of $F$. In the similar way, the covariant functor  $\hm_{\mathscr{C}}(-,G(-)):\mathscr{C}^{op}\times\mathscr{D}\rightarrow\mathbf{Set}$ is called the right hom-set adjunction of $G$.  \\

\begin{definition} Let $F:\mathscr{C}\rightarrow\mathscr{D}$  and $G:\mathscr{D}\rightarrow\mathscr{C}$ are two covariant functors. The functor $F$ is called a left adjoint to the functor $G$ (or $G$ is called a right adjoint to $F$) If there exists a natural isomorphism $\Phi:\hm_{\mathscr{D}}(F(-),-)\Rightarrow\hm_{\mathscr{C}}(-,G(-))$ between the left and right hom-set adjunctions of $F$ and $G$ respectively. Sometimes, this relationship is indicated by $\Phi:F\dashv G$.
\end{definition}

\begin{definition}  Let $\mathscr{C}$ be a  subcategory of the category $\mathscr{D}$. Then $\mathscr{C}$ is called a reflective subcategory of $\mathscr{D}$, if there exists a covariant functor $F:\mathscr{D}\rightarrow\mathscr{C}$  which is a left adjoint to the inclusion functor
$i:\mathscr{C}\rightarrow\mathscr{D}$, $i.e.,$ $F\dashv i$.  In this situation, the functor $F$ is called a reflector.
\end{definition}

\section{Main Results}
Although the following result is probably known, we have provided a proof for the sake of completeness. \\

\begin{proposition} \label{correspond} For any ring $R$, consider the space $X=\Spec(R)$ with the Zariski topology. Then there exists a one-one correspondence between $I(R)$ and $\Cl(X)$.
 \end{proposition}
{\bf Proof.} Define the map $\eta:I(R)\rightarrow\Cl(X)$ by $\eta(e)= D(e)$ for any $e\in I(R)$. Since  $D(e)=V(1-e)$, hence $D(e)$ is a clopen and so $\eta$ is well-defined. We shall show that it is bijective. For the injectivity, let that $\eta(e)=\eta(e')$ for some $e,e'\in I(R) $, then $V(1-e)=V(1-e')$. Hence, $\sqrt{\langle 1-e\rangle}=\sqrt{\langle 1-e'\rangle}$, so $1-e=(1-e)^{n}\in \langle 1-e'\rangle$ for some $n\geq 1$. Thus $1-e=r'(1-e')$ for some $r'\in R$,  if we multiply this equation to $e'$ we get $e'=ee'$. By the similar way we also get $e=ee'$ thus $e=e'$.\\
Now we show that $\eta$ is surjective. Let that $U$ be a clopen of $X$. Since $U$
is closed it is quasi-compact   and similarly  its complement. Write $U=\bigcup\limits_{i=1}^{n}D(f_{i})$ as a finite union of standard opens. Similarly, write $\Spec(R)\setminus U=\bigcup\limits_{j=1}^{m}D(g_{j})$ as a finite union of standard opens. But $D(f_{i}g_{j})= D(f_{i})\cap D(g_{j})=\emptyset$, and so each of $f_{i}g_{j}$ is nilpotent, thus if we set $I=\langle f_{1},...,f_{n}\rangle$ , $J=\langle g_{1},...,g_{m}\rangle$ then $(IJ)^{N}=0$ for some  sufficiently large integer $N\geq 1$. Also $\Spec(R)=V(I^{N})\coprod V(J^{N})$  because $V(I^{N})=V(I)=\bigcup\limits_{j=1}^{m}D(g_{j})$ and similarly $V(J^{N})=V(J)=\bigcup\limits_{i=1}^{n}D(f_{i})$  so we have $R=I^{N}+J^{N}$. Write $1=x+y$ with $x\in I^{N}$ and $y\in J^{N}$. Then we have $x=x^{2}$ because $x-x^{2}=x(1-x)=xy=0$ for the last equality note that $xy\in I^{N}J^{N}=(IJ)^{N}=0$. So $x$ is an idempotent and we have $\eta(x)=D(x)=\bigcup\limits_{i=1}^{n}D(f_{i})=U$. $\Box$  \\

For a ring $R$, the elements $0$ and $1$ are called the trivial idempotents.  \\

\begin{corollary} \label{idemconnected} For any ring $R$, set $X=\Spec(R)$ with the Zariski topology. Then the space $X$ is connected if and only if the idempotents of $R$ are only $0$ and $1$.
\end{corollary}
{\bf Proof.} This is a direct consequence of the above proposition.  $\Box$ \\

\begin{proposition} \label{goodlem} Let $M$ be a regular ideal of $R$. Then $M$ is a max-regular ideal if and only if the idempotents of $R/M$ are trivial.
 \end{proposition}
{\bf Proof.} If the set of idempotents of $R/M$ is trivial then it is easy to see that $M$ is a max-regular ideal of $R$.\\ Conversely, suppose that $M$ be a max-regular ideal of $R$. Let $\overline{x}=x+M$ be an arbitrary idempotent of $R/M$ where $x \in R$, this implies that $x-x^{2}\in M$. By the lemma \ref{regularidem}, there exists an idempotent $e \in M$ so that $(x-x^{2})=e(x-x^{2})$ thus we get $(1-e)(x-x^{2})=0$. This implies that $(1-e)x$ is an idempotent of $R$. Write $x= (1-e)x+ ex$ which belongs to the regular ideal $M+\langle(1-e)x\rangle$. Also
 \begin{equation}\label{equa} M\subseteq M+\langle(1-e)x\rangle\subseteq M+\langle x\rangle. \end{equation}  We have $M+\langle x\rangle=R$ or $M+\langle x\rangle\neq R$. \\

 if $M+\langle x\rangle=R$, then write $1=rx+r_{1}e_{1}+...+r_{n}e_{n}$ where $r,r_{i}\in R$ and  $e_{i}\in M$, if we multiply this equation to $1-x$ then we get $1-x=r(x-x^{2})+(r_{1}e_{1}+...+r_{n}e_{n})(1-x)$ which belongs to $M$ and so $\overline{x}=1$, which is a trivial idempotent. \\

 If $M+\langle x\rangle\neq R$, then since $M$ is a max-regular ideal and $ M+\langle(1-e)x\rangle$ is a regular ideal so from the (\ref{equa}) we get $M= M+\langle(1-e)x\rangle$ but since $x \in  M+\langle(1-e)x\rangle=M $ so $\overline{x}=0$, therefore in this case also $\overline{x}$ is a trivial idempotent.  Therefore the idempotents of $R/M$ are only trivial. $\Box$ \\

\begin{corollary} \label{above}Let that $M$ be a max-regular ideal of $R$. Then $V(M)$ is a connected subset of $X=\Spec(R)$.
\end{corollary}
{\bf Proof.} By the proposition \ref{goodlem}, the idempotents of $R/M$ are trivial and so by the corollary \ref{idemconnected}, the space $\Spec(R/M)$ is connected. In other hand, $V(M)$ is naturally homeomorphic to the $\Spec(R/M)$, so it is also connected. $\Box$

\begin{theorem}\label{max-reg} Let $R$ be a ring and set $X=\Spec(R)$ with the Zariski topology. Then $C\subseteq X$ is
a connected component if and only if  $C$ is of the form $V(M)$ for some max-regular ideal $M$ of $R$.
\end{theorem}

{\bf Proof.} In order to prove the assertion, first we define a map \\ $f: X\rightarrow \mr(R)$ by $f(\frak{p})=\langle e$ $|$ $e\in \mathfrak{p}\cap I(R)\rangle$. It is easy to check that for any prime $\mathfrak{p}$,  $f(\mathfrak{p})$ is a max-regular ideal of $R$ and so the map is well defined. Also $f$ is continuous, because for any basis open  $\mathscr{O}_{e}$ of $\mr(R)$ where $e\in I(R)$, we have $f^{-1}(\mathscr{O}_{e})=D(e)$.  \\

 Now let that $C$ is a connected component of $X$ then $f(C)$ is connected subset of $\Sp(R)$, because $f$ is continuous. But by the lemma \ref{mrprofinite}, $\mr(R)$ is a profinite space and so $f(C)=\{M\}$ for a max-regular ideal $M$ of $R$. But we have  $C\subseteq f^{-1}(\{M\})=V(M)$.  Also by the  corollary \ref{above}, $V(M)$ is a connected subset of $X$ so the inclusion $C\subseteq V(M)$ implies the equality, because $C$ is a connected component.\\

Conversely, assume that $M$ be a max-regular ideal of $R$.  Again by the corollary \ref{above}, $V(M)$ is a connected subset of $X$, thus it is contained in a connected component $C$  of $X$. By the above paragraph, $C= V(N)$ for some max-regular ideal $N$ of $R$. So $V(M)\subseteq V(N)$ this implies that $N\subseteq\sqrt{N}\subseteq \sqrt{M}$. But for each  element $e$ of a set of  idempotent generators of $N$, we have $e=e^{n}\in M$ for some $n\geq 1$. Hence,
$N\subseteq M$. Since $N$ is a max-regular ideal, so $N=M$.   $\Box$  \\

\begin{corollary}\label{coarser} For a ring $R$, set $X=\Spec(R)$ with the Zariski topology. Then there exists a topology $\mathscr{T}$ on $X/_{\sim}$ which is coarser than the quotient topology and the space $(X/_{\sim}, \mathscr{T})$ is  profinite.
\end{corollary}
{\bf Proof.}  In the light of the theorem \ref{max-reg}, we have   $X/_{\sim}=\{V(M)$ $|$ $M\in\mr(R) \}.$
 Hence the map $\Phi:X/_{\sim}\rightarrow\mr(R)$ given by  $V(M)\rightsquigarrow M$ is well-defined and bijective. Furthermore, by the remark \ref{fineremark}, there exists a basis for $\mr(R)$ in which any element  of this basis is of the form $\mathscr{O}_{e}=\{M\in\mr(R)$ $|$ $e\notin M\}$  where $e\in I(R)$. So the map $\Phi$ induces a topology $\mathscr{T}$ on $X/_{\sim}$ with the basis $\{\Phi^{-1}(\mathscr{O}_{e})$ $|$ $e\in I(R)\}$. Therefore with this topology, $\Phi$ is a homeomorphism. Also the topology $\mathscr{T}$ is coarser than  the quotient topology, because $\pi^{-1}(\Phi^{-1}(\mathscr{O}_{e}))= D(e)$ where $\pi:X\rightarrow X/_{\sim} $ is the canonical projection. By the lemma \ref{mrprofinite}, the space $\mr(R)$  is  profinite,  hence the space $(X/_{\sim},\mathscr{T})$ is also profinite.    $\Box$ \\

For given ring $R$, denote by  $\mathfrak{M}(R)$ the set of all maximal ideals of $R$ and consider it as a topological subspace of  $\Spec(R)$.  \\

\begin{theorem} Let $X$ be a topological space and set $R=C(X)$ the ring of real-valued continuous functions on $X$. Then there exists a topology $\mathscr{T}$  on the set  $\mathfrak{M}(R)/_{\sim}$ which is coarser than the quotient topology and the space $(\mathfrak{M}(R)/_{\sim},\mathscr{T})$ is profinite.
\end{theorem}

{\bf Proof.}  First we show that the map $\Psi:\Spec(R)/_{\sim}\rightarrow\mathfrak{M}(R)/_{\sim}$ given by $V(M)\rightsquigarrow V(M)\cap\mathfrak{M}(R)$ is well-defined and bijective. \\

To achieve the purpose we act as follows, according to \cite[Theorem  2.11 ]{Gillman}, every prime ideal $\mathfrak{p}$ of $R$ is contained in a unique maximal ideal $\mathfrak{m}_{\mathfrak{p}}$. Hence, we obtain a map $\psi:\Spec(R)\rightarrow\mathfrak{M}(R)$ given by $\psi(\mathfrak{p})=\mathfrak{m}_{\mathfrak{p}}$. This map is continuous according to \cite[Corollary 1.6.2.1]{Bkouche}. Therefore, $\psi$ maps any connected component $V(M)$ of $\Spec(R)$ to a connected subset $\psi(V(M))=V(M)\cap\mathfrak{M}(R)$ of $\mathfrak{M}(R)$. In fact, we show that $V(M)\cap\mathfrak{M}(R)$ is  a connected component of $\mathfrak{M}(R)$. For, choose
 $\mathfrak{m}\in V(M)\cap\mathfrak{M}(R)$ and let $C$ is a connected component of $\mathfrak{M}(R)$ containing $\mathfrak{m}$, since the inclusion map $\mathfrak{M}(R)\hookrightarrow\Spec(R)$ is continuous then $C$
is a connected subset of $\Spec(R)$, so $C$ is contained in a connected component of $\Spec(R)$ which is exactly $V(M)$, because $\mathfrak{m}\in V(M)$. Thus, we have $C\subseteq V(M)\cap\mathfrak{M}(R)$, but by the connectedness of $V(M)\cap\mathfrak{M}(R)$ we obtain $C= V(M)\cap\mathfrak{M}(R)$. \\

Therefore, the map $\Psi:\Spec(R)/_{\sim}\rightarrow\mathfrak{M}(R)/_{\sim}$ is well-defined (note that the map $\Psi$ is in fact  induced by $\psi$). Surjectivity of $\Psi$ is clear from the preceding argument.
For its injectivity, suppose that $V(M)\cap\mathfrak{M}(R)=V(N)\cap\mathfrak{M}(R)$ for some max-regular ideals $M$ and $N$ of $R$. If $M\neq N$ then we can choose an idempotent $e\in M\setminus N$. Let $\mathfrak{m}$ be a maximal ideal of $R$ containing $M$, then $\mathfrak{m}$ also contains $N$ and so the regular ideal $N+\langle e\rangle$ is contained in $\mathfrak{m}$, but since $N$ is a max-regular ideal we get $N=N+\langle e\rangle$ which is a contradiction. \\

  Finally, by using the corollary \ref{coarser}, the bijective map  $\Psi:\Spec(R)/_{\sim}\rightarrow\mathfrak{M}(R)/_{\sim}$ induces a topology $\{\Psi(V)$ $|$ $V\in\mathscr{T}\}$ ( we denote it also by $\mathscr{T}$) on $\mathfrak{M}(R)/_{\sim}$, and with this topology the map $\Psi$ becomes   a homeomorphism and so the space $(\mathfrak{M}(R)/_{\sim},\mathscr{T})$ is profinite.  $\Box$  \\

\begin{theorem} \label{cfunction} If  $X$ is a compact Hausdorff space, then the set $X/_{\sim}$ endowed with some topology $\mathscr{T}$ is a profinite space.

\end{theorem}

 {\bf Proof.} For a compact  Hausdorff space $X$, according to \cite[4.9.(a)]{Gillman}, the map $\mu:X\rightarrow\mathfrak{M}(R)$ given by  $\mu(x)=\mathfrak{m}_{x}=\{f\in R$ $|$ $f(x)=0\}$ is a homeomorphism where $R=C(X)$. Finally, the result implies from the preceding theorem.  $\Box$  \\

\begin{remark}\label{Dehghan} Let $g:X\rightarrow Y$ be a continuous map where $X$ is compact Hausdorff and $Y$ is Hausdorff. Then $g$ is a closed map, because each closed subset $F$ of $X$ is compact and so $g(F)$ is a compact subset in $Y$, but Hausdorffness of $Y$ implies that $g(F)$ is closed.  \\
\end{remark}

\begin{remark} For given topological space $X$, set $R=C(X)$ the ring of real-valued continuous functions. Then one can easily check that the set of idempotents $I(R)$ of $R$ is exactly equal to the set  $\{\chi_{_{U}} : U\in\Cl(X) \}$ where $\chi_{_{U}}$ is the characteristic function of $U$. In the sequel we need to this characterization of the idempotents. Also this characterization of idempotents implies that $X$ is connected if and only if the space $\Spec(R)$ is connected.  \\
\end{remark}

  The above theorem \ref{cfunction}, leads us to a covariant functor $F:\mathscr{K}\rightarrow\mathscr{P}$ from the category of compact  Hausdorff spaces $\mathscr{K}$ to the category of profinite spaces $\mathscr{P}$ and this categorical construction  implies that  the category of profinite spaces is a reflective subcategory of the category of compact  Hausdorff spaces and the reflector is the foregoing functor.  Hence in what follows,  we plan to describe this  functor explicitly and then show that this functor actually is a reflector.\\

 \begin{remark}\label{functor} For any compact  Hausdorff space $X$, set $R=C(X)$ the ring of real-valued continuous functions on $X$, also set $\mathfrak{m}_{x}=\{f\in R$ $|$ $f(x)=0\}$ the maximal ideal of $R$ corresponding to each $x\in X$. Then by the theorem \ref{cfunction}, each connected component of  $X$ is of the form $C_{M}=\{x\in X$ $|$ $M\subseteq\mathfrak{m}_{x}\}$ where $M\in\mr(R)$. \\

 Therefore, $X/_{\sim}=\{C_{M}$ $|$ $M\in \mr(R) \}$ and the subsets  $\mathscr{O}_{U}=\{C_{M}$ $|$ $\chi_{_{U}}\notin M\}$ where $\chi_{_{U}}$ is the characteristic function of  $U\in \Cl(X)$, as open subsets constitute a basis for the  topology $\mathscr{T}$ as given in the theorem \ref{cfunction}.\\

Finally, define the functor $F:\mathscr{K}\rightarrow\mathscr{P}$,  for each compact Hausdorff space $X$, by $F(X)=(X/_{\sim},\mathscr{T})$.
 Moreover, for any continuous function $f:X\rightarrow Y$ between the compact  Hausdorff spaces, then $Ff: (X/_{\sim},\mathscr{T})\rightarrow (Y/_{\sim},\mathscr{T'})$ is defined for each $C_{M}\in X/_{\sim}$, by $(Ff)(C_{M})= C_{N}$ where  $C_{N}\in Y/_{\sim}$ and $f(C_{M})\subseteq C_{N}$. Consider  the following commutative diagram
 $$\xymatrix{
X \ar[r]^{f} \ar[d]^{\pi_{X}} & Y \ar[d]^{\pi_{Y}} \\ (X/_{\sim},\mathscr{T}) \ar[r]^{Ff} & (Y/_{\sim},\mathscr{T'}) }$$

where the vertical arrows are the canonical projections. The following lemma guarantees that $F$ actually is a functor.  \\
\end{remark}

\begin{lemma} With the notation and assumptions as in the previous remark, then the map $Ff$ is continuous.
\end{lemma}

 {\bf Proof.} For continuity of $Ff$ it is enough to show that $(Ff)^{-1}(\mathscr{O}_{V})=\mathscr{O}_{f^{-1}(V)}$ for each open basis $\mathscr{O}_{V}$ of $\mathscr{T'}$. \\

First, suppose that $C_{M}\in \mathscr{O}_{f^{-1}(V)}$, then we have $\chi_{_{f^{-1}(V)}}\notin M$. Also, $(Ff)(C_{M})=C_{N}$ where $f(C_{M})\subseteq C_{N}$.
But the max-regular ideal $M\subseteq C(X)$ is contained in a maximal ideal $\mathfrak{m}_{x}$ of $C(X)$ for some $x\in X$. Thus $x\in C_{M}$ so $y=f(x)\in f(C_{M})\subseteq C_{N}$.  Hence, $N\subseteq\mathfrak{m}_{y}$ where $\mathfrak{m}_{y}=(f^{\ast})^{-1}(\mathfrak{m}_{x})$ for which the ring map $f^{\ast}:C(Y)\rightarrow C(X)$ is induced by $f:X\rightarrow Y$. But from the max-regular property of $M$ we get $f^{\ast}(\chi_{_{V}})=\chi_{_{f^{-1}(V)}}\notin \mathfrak{m}_{x}$. Therefore, $\chi_{_{V}}\notin\mathfrak{m}_{y}$, hence $\chi_{_{V}}\notin N$ and so $C_{N}\in\mathscr{O}_{V} $. \\

For the inverse inclusion, suppose that $C_{M}\in (Ff)^{-1}(\mathscr{O}_{V})$. Then $(Ff)(C_{M})=C_{N}\in \mathscr{O}_{V}$ where $f(C_{M})\subseteq C_{N}$. Thus $\chi_{_{V}}\notin N$. As above, suppose that $M\subseteq \mathfrak{m}_{x}$ for some maximal ideal $\mathfrak{m}_{x}$ of $C(X)$. Hence $x\in C_{M}$ and so  $y=f(x)\in f(C_{M})\subseteq C_{N}$. Thus $N\subseteq\mathfrak{m}_{y}=(f^{\ast})^{-1}(\mathfrak{m}_{x})$. So, by the max-regular property of $N\subseteq C(Y)$ we have $\chi_{_{V}}\notin\mathfrak{m}_{y}$. Thus $f^{\ast}(\chi_{_{V}})=\chi_{_{f^{-1}(V)}}\notin\mathfrak{m}_{x}$. Therefore $\chi_{_{f^{-1}(V)}}\notin M$ and so $C_{M}\in\mathscr{O}_{f^{-1}(V)}$. $\Box$  \\

\begin{theorem} \label{alternative} The category of profinite spaces is a reflective subcategory in the category of compact Hausdorff spaces.
\end{theorem}

{\bf Proof.}  We prove that the functor $F:\mathscr{K}\rightarrow\mathscr{P}$ as defined in the remark \ref{functor}, is a left adjoint to the inclusion functor
$i:\mathscr{P}\rightarrow\mathscr{K}$. For this purpose, we show that there exists a natural isomorphism $\mu$ between the left and right hom-set adjunctions of $F$ and $i$ respectively, \\ $\mu:\hm(F(-),-)\Rightarrow\hm(-,i(-))$ . \\

 For each object $(X,P)\in \mathscr{K}^{op}\times\mathscr{P} $ we define the natural transformation
$\mu_{_{X,P}}:\hm((X/_{\sim},\mathscr{T}),P)\rightarrow\hm(X,P)$
by $\mu_{_{X,P}}(g)=g\circ\pi_{X}$ where $g\in\hm((X/_{\sim},\mathscr{T}),P)$ and  $\pi_{X}:X\rightarrow X/_{\sim}$ is the canonical projection. The map $\mu_{_{X,P}}$ is injective because $\pi_{X}$ is surjective.\\

For the surjectivity of $\mu_{_{X,P}}$, suppose that $f\in\hm(X,P)$.
Put $g=\pi_{P}^{-1} \circ Ff :(X/_{\sim},\mathscr{T})\rightarrow P$ (note that if $P$ is already profinite then the natural projection $\pi_{P}: P\rightarrow (P/_{\sim},\mathscr{T'})$ by the remark \ref{Dehghan}, is a homeomorphism). Finally, commutativity of the following diagram\\

$$\xymatrix{
X \ar[r]^{f} \ar[d]^{\pi_{X}} & P \ar[d]^{\pi_{P}} \\ (X/_{\sim},\mathscr{T}) \ar[r]^{Ff} & (P/_{\sim},\mathscr{T'}) }$$

implies that $\mu_{_{X,P}}(g)=g\circ\pi_{X}=f$.   $\Box$ \\

\textbf{On the connected components of t(X).} Denote by \textbf{Top} the category of topological spaces with continuous maps as morphisms. In what  follows, we recall the definition of the  classical covariant functor $t:\textbf{Top} \rightarrow \textbf{Top}$ and then state  some of its  basic properties.

\begin{definition}\label{tfunctor} The functor $t:\textbf{Top} \rightarrow \textbf{Top}$ is defined for a topological space $X$, by $t(X)=\{Z\subseteq X$ $|$  $Z$ is a closed and irreducible subset of $X\}$. \\ The  subsets  $t(Y)$ where $Y$ is a closed subset of $X$, as closed subsets constitute a topology for $t(X)$. \\ Moreover, for any continuous map $f:X\rightarrow X^{'}$, then the map $tf:t(X)\rightarrow t(X^{'})$ is defined for each $Z\in t(X)$ by $(tf)(Z)=\overline{f(Z)}$  which is a continuous map.
\end{definition}

We will need the following easy lemma.  \\

\begin{lemma}\label{imconnect} Let $X$ be any topological space and let $C$ be a connected component of it, then $t(C)$ is a connected subset of $t(X)$.
 \end{lemma}
{\bf Proof.} Suppose that $t(C)=(\mathcal{U}\cap t(C))\cup (\mathcal{V} \cap t(C))$ be a disjoint separation for $t(C)$ where $\mathcal{U}$ and $\mathcal{V}$ are open subsets of $t(X)$. Set $\mathcal{U}=t(X)\setminus t(E)$ and $\mathcal{V}=t(X)\setminus t(F)$ where $E$ and $F$ are closed subsets of $X$. Also Set $U=X \setminus E$ and $V=X \setminus F$, then $C=(U \cap C) \cup (V \cap C)$  is a disjoint separation by the open subsets of $C$. But by the connectedness of $C$ we get $C=U\cap C$ or $C=V \cap C$. Thus we get  $t(C)= \mathcal{U} \cap t(C)$ or $t(C)= \mathcal{V} \cap t(C)$.     $\Box$

\begin{remark} Note that in the above lemma the assertion is also true for any connected subset of $X$. More precisely, for any connected subset $C$ of
$X$ then the set $\{Z\in t(X)$ $|$ $Z\subseteq C\}$ is a connected subset of $t(X)$. The proof is similar to the proof of the above
 lemma.

\end{remark}

\begin{remark}  Structure of the connected components of $t(X)$ in the general case is as follows. Each topological space $X$ can be written as a disjoint union of the its connected components, i.e., $X=\coprod\limits_{i\in I} C_{i}$ where  each $C_{i}$ is a connected component of $X$. From this fact we easily get   $t(X)=\coprod\limits_{i\in I} t(C_{i})$. Hence,  each connected component $\mathscr{C}$ of $t(X)$,  is of the form
$\mathscr{C}=\coprod\limits_{j\in J} t(C_{j})$ for some $J\subseteq I$, because by the lemma \ref{imconnect}, each of the $t(C_{i})$ is connected.\\
 Now if $J$ is a finite set then it is just a single point subset. Namely, $\mathscr{C}=t(C_{j})$ where $J=\{j\}$, because in the finite case each of the  $t(C_{j})$ is a disjoint open subset of $\mathscr{C}$.\\ However, in the general case the set $J$ is not necessarily finite. \\ The following theorem says us that the set $J$ is finite whenever  the space $X/_{\sim}$ is totally disconnected with some topology $\mathscr{T}$.  \\
\end{remark}

\begin{theorem} \label{conn component}Let $X$ be a topological space so that $(X/_{\sim},\mathscr{T})$ is totally disconnected  with some topology $\mathscr{T}$ which is coarser than the quotient topology. Then $\mathscr{C}\subseteq t(X)$ is a connected component of  $t(X)$ if and only if $\mathscr{C}$ is of the form $t(C)$ where $C$ is a connected component of $X$.
 \end{theorem}
 {\bf Proof.}  First we should note that any closed and irreducible subset $Z$ of $X$ is also a connected subset of $X$. We denote by $\Gamma(Z)$ the connected component of $X$ which contains $Z$. So this defines a map \\  $\Gamma:t(X)\rightarrow (X/_{\sim},\mathscr{T})$. The map $\Gamma$ is  continuous. In order to prove this, for any closed subset $\mathscr{E}$  of $X/_{\sim}$, we will show that $\Gamma^{-1}(\mathscr{E})=t(E)$ where $E= \pi^{-1}(\mathscr{E})$ and $\pi:X\rightarrow X/_{\sim}$ is the canonical projection.  \\

 First let that $Z \in t(E)$, since $Z$ is a nonempty connected subset thus
$Z\subseteq \pi(z_{0})=\Gamma(Z)$ for some point $z_{0} \in Z$, but $\pi(z_{0}) \in \mathscr{E}$ because $Z\subseteq E$ therefore $\Gamma(Z) \in \mathscr{E} $ this shows that $t(E)\subseteq \Gamma^{-1}(\mathscr{E})$. For the inverse inclusion, suppose that $Z \in \Gamma^{-1}(\mathscr{E})$ then $\Gamma(Z) \in \mathscr{E}$, but $Z$ is connected, so for any point $z \in Z$ we have $\pi(z)= \Gamma(Z) \in \mathscr{E}$ this shows that $z \in E=\pi^{-1}(\mathscr{E})$ thus $Z\subseteq E$ this means that $Z \in t(E)$.\\

Now for proving the assertion, let that $\mathscr{C}$ be any connected component of $t(X)$. Since the map $\Gamma$ is continuous and the space $(X/_{\sim},\mathscr{T})$ is totally disconnected, then  we have $\Gamma(\mathscr{C})=\{C\}$ for some single point subset $\{C\}$ of $X/_{\sim}$ where $C$ is a connected component of $X$. But $\mathscr{C} \subseteq t(C)$, because for any point $Z \in \mathscr{C}$ we have $\Gamma(Z)=C$, thus this means that $Z \in t(C)$. But by the lemma \ref{imconnect},  $t(C)$ is a connected subset of $t(X)$, also $\mathscr{C}$ is a connected component of $t(X)$, thus we get $\mathscr{C}=t(C)$.\\
 Conversely, assume that $C$ is a connected component of $X$. By the lemma \ref{imconnect}, $t(C)$
 is a connected subset of $t(X)$, thus $t(C)$ is contained in some connected component $\mathscr{C}$ of $t(X)$. But we have $\{C\}=\Gamma(t(C)) \subseteq \Gamma(\mathscr{C})=\{C^{'}\}$ for some $C^{'} \in X/_{\sim}$. Thus we get $C=C^{'}$ and $t(C)=\mathscr{C}$, so $t(C)$ is a connected component.    $\Box$ \\

\begin{definition}\label{sober} Let $X$ be a topological space. If any closed and irreducible subset $Z$ of $X$ have a unique generic point $\zeta$ ( i.e., there exists only a unique point $\zeta\in Z$ with the property $Z=\overline{\{\zeta\}}$ ) then  $X$ is called a sober space.\\
\end{definition}

For a ring $R$, the set $X=\Spec(R)$ with  the Zariski topology is a sober space. Because, any closed and irreducible subset of $X$ is of the form $V(\mathfrak{p})=\overline{\{\mathfrak{p}\}}$. As another example, totally disconnected spaces are sober, because a closed and irreducible subset of any topological space is connected. The following proposition will provide more examples.\\

\begin{proposition} For any topological space $X$ then we have,\\
\begin{itemize}
  \item[(i)] There exists a natural transformation of functors $\alpha:id\Rightarrow t$ where $id:\textbf{Top} \rightarrow \textbf{Top}$ is the identity functor and for the functor $t$ see the definition $\ref{tfunctor}$.

 \item[(ii)] The rule $E\rightsquigarrow t(E)$ is a bijective between  the closed subsets of $X$ and the closed subsets of $t(X)$.
\item[(iii)] The space $t(X)$ is sober.

\end{itemize}
\end{proposition}
 {\bf Proof.}\begin{itemize}  \item[(i)] For each object $X\in\textbf{Top}$, define the continuous map $\alpha_{X}:X\rightarrow t(X)$ given by  $x\rightsquigarrow\overline{\{x\}}$. For a continuous map $f:X\rightarrow X'$ we have $\overline{\{f(x)\}}=\overline{f(\overline{\{x\}})}$, which implies that the following diagram is commutative.

                  $$\xymatrix{
X \ar[r]^{\alpha_{X}} \ar[d] & t(X) \ar[d] \\ X' \ar[r]^{\alpha_{X'}} & t(X') }$$

and so $\alpha:id\Rightarrow t$ is a natural transformation of functors. \\
 \
 \item[(ii)] Injectivity of the map is easy and its surjectivity implies from the definition \ref{tfunctor}. \\

 \item[(iii)] Any closed and irreducible subset of  $t(X)$ is of the form $t(Z)$ where $Z$ is a closed and irreducible subset of $X$. Also $t(Z)=\overline{\{Z\}}$ where the closer is taken in $t(X)$. Now if there exists another one $Z' \in t(Z)$ for which
     $t(Z)=\overline{\{Z'\}}$ then $t(Z)=t(Z')$. But by the part (ii) above, one has $Z=Z'$. $\Box$  \\

 \end{itemize}

\end{document}